\documentclass[clean,reqno,twoside,11pt]{amsart}
\usepackage{mltools}
\usepackage{math62,mlthm10}

\usepackage{amscd}
\usepackage{amsmath} 
\usepackage{amssymb,amsfonts,mathrsfs}

\usepackage[mathscr]{eucal}
\usepackage{bm}

\def\mlSetRCSInfo#1{\def\mlRCSInfo{#1}}
 
\mlSetRCSInfo{\today}
 
\newcommand{\T}{\mathbb{T}}

\def\lapa{\triangle}

\def\perta{\triangle_k}
\def\pertb{\triangle^{(0, 1)}_k}

\usepackage{LesMos_local}
\allowdisplaybreaks
\begin{document}

\title{Modular Gaussian curvature} 

\author{Matthias Lesch}
\address{Mathematisches Institut,
Universit\"at Bonn,
Endenicher Allee 60,
53115 Bonn,
Germany}

\email{ml@matthiaslesch.de, lesch@math.uni-bonn.de}
\urladdr{www.matthiaslesch.de, www.math.uni-bonn.de/people/lesch}
\thanks{The work of M.L. was partially supported by the 
        Hausdorff Center for Mathematics, Bonn}

\author{Henri Moscovici}
\address{Department of Mathematics,
The Ohio State University,
Columbus, OH 43210,
USA}
\email{henri@math.ohio-state.edu}

\thanks{The work of H. M. was partially
 supported by the National Science Foundation
    award DMS-1600541} 
\dedicatory{Dedicated to Alain Connes with admiration and much appreciation}
 
\maketitle 
 
\date{\today}


\section{Introduction}

The genesis of natural but noncommuting coordinates can be traced back to
 Heisenberg's uncertainty principle in quantum mechanics, which limits the
accuracy of the simultaneous determination of the position and momentum $(q,p)$
of a subatomic particle. As Heisenberg argued~\cite{Hei1927} (and
Kennard rigorously derived~\cite{Ken1927}), the inherent  
imprecision of such a measurement ``is a straightforward mathematical consequence
of the quantum mechanical commutation rule for the position and the corresponding
momentum operators $q p - p q = i \hbar$'', 
where $ \hbar= \frac{h}{2\pi}$ is the reduced Planck constant.
Such an identity cannot be satisfied by matrices (over $\C$), which is obvious,
but not even by bounded operators in Hilbert space. Assuming $q$ and $p$ self-adjoint, 
this can be seen by passing to the Weyl integrated form~\cite[\S 45]{Wey28}, 
\begin{align} \label{H-def}
V_s U_t = e^{2 \pi i \hbar ts } U_t V_s,  \qquad t, s \in \R. 
\end{align}
Moreover,
the latter relation determines a unitary representation $\pi_{\hbar}$ of the (implicitly defined)
group $H_3 (\R)$, called the Heisenberg group. By a celebrated theorem of
Stone and von Neumann all such irreducible representations are unitarily equivalent. 
The restriction to the lattice $H_3 (\Z) \subset H_3 (\R)$ of an irreducible unitary representation 
$\pi_{\theta}$, $\theta \in \R$, generates the $C^*$-algebra $A_\theta$ 
nowadays known as the  {\em  noncommutative torus} of slope $\theta$. 
When $\theta \in [0,1] \setminus \Q$, as shall be assumed throughout this paper,
the $C^*$-algebra $A_\theta$ is also known as the  {\em irrational rotation algebra} and is  
the (unique up to isomorphism) $C^*$-algebra generated by a pair of unitary operators 
$\, U_1, \,  U_2 $ satisfying
\begin{align} \label{A-def}
U_2U_1=e^{2 \pi i \theta} U_1U_2 .
\end{align}
Moreover $A_\theta$ is a simple $C^*$-algebra, thus
typifying the coordinates of a purely noncommutative space.
For this reason on the one hand, and due to its accessibility on the
other, $A_\theta$ has received much attention during the
last several decades and has been a favorite testing ground for quite a number of
fruitful mathematical investigations. 

Although the habitual geometric intuition is rendered utterly
inoperative in a `space without points' such as the one represented by $A_\theta$, the
curvature as a ``measure of deviation from flatness'' (in Riemann's
own words) could still make some sense.  It is the goal of this brief
survey to review the recent developments that led to the emergence of
a quantized version of Gaussian curvature for the
noncommutative torus.  
Many of the essential ideas presented below
have their origin in Alain Connes' 1980 C.~R.~Acad.~Sc.~Paris
Note~\cite{Con:CAG}, which effectively constitutes the birth
certificate of noncommutative differential geometry. That foundational
article not only established the most basic geometric concepts and constructions,
such as the geometric realization of the finitely generated projective modules over
$A_\theta$, the explicit construction of constant curvature
connections for them and the definition and calculation of
their Chern classes, but also provided  the crucial computational tool, 
in the form of a pseudodifferential calculus adapted to $C^*$-dynamical systems.
  
The specific line of research whose highlights we are about to
summarize was sparked by a paper by Connes and Paula
Cohen ({\it Conformal geometry of the
  irrational rotation algebra}, Preprint MPI Bonn, 1992-93)
 which showed how the passage
from the (unique) trace of $A_\theta$ to a non-tracial conformal
weight associated to a Weyl factor (or `dilaton') gives rise to a non-flat geometry on
the noncommutative torus, which can be investigated with the
help of the adapted pseudodifferential calculus of~\cite{Con:CAG}. In
a later elaboration~\cite{ConTre2011} of that paper, the passage from flatness to
conformal flatness was placed in the setting of spectral triples 
(see \S \ref{flatST} and \S \ref{ModST} below),
which in the intervening years has emerged as the proper framework for
the metric aspect in noncommutative geometry (cf.~\cite[Ch.~6]{Con:NG},
\cite{Con2013}, \cite{cmIII}). Completing the calculations begun in the 1992 preprint
they proved in~\cite{ConTre2011} an analogue of the Gauss-Bonnet
formula for the conformally twisted (called `modular') spectral triples. The
full calculation of the modular Gaussian curvature was first done by
A. Connes in 2009, with the aid of Wolfram's Mathematica, and is
included in~\cite{CM2014}. Fathizadeh and Khalkhali~\cite{FarKha2013} independently performed the same calculation with the help of a different computing software.  

Apart from computing the expression of the modular curvature (see \S \ref{ModC} below),
Connes and Moscovici showed in \cite{CM2014} 
that one can make effective use of variational
methods even in the abstract operator-theoretic context of the spectral triple encoding
the geometry of the noncommutative torus.  After giving a 
variational proof of the modular Gauss-Bonnet formula 
 which requires no computations (see \S \ref{ModGB}),
they related the modular Gaussian curvature to the gradient of the Ray-Singer
log-determinant of the Laplacian viewed as a functional on the space
of Weyl factors. 
As a consequence, they obtained an a priori proof
of an internal consistency relation for the constituents of the
modular curvature. In addition they showed by purely
operator-theoretic arguments that, as in the case of Riemann surfaces
(cf.~\cite{OPS88}), the normalized log-determinant functional attains
its extreme value only at the trivial Weyl factor, in other words for the flat `metric'
 (see \S \ref{var-det} below). 
 
For reasons which will soon become transparent (see \S \ref{FolAlg}),
the natural equivalence relation between noncommutative spaces  is that
of Morita equivalence between their respective algebras of coordinates. 
For noncommutative tori the Morita equivalence is implemented by the
Heisenberg bimodules described by Connes~\cite{Con:CAG} and Rieffel~\cite{Rie1981}.
Lesch and Moscovici extended in~\cite{LM2015} the  
results of \cite{CM2014} to spectral triples on noncommutative tori associated to 
Heisenberg equivalence bimodules (see \S \ref{ModHST} and \S \ref{RSYM}). Moreover,
in doing so they managed to dispose of any computer-aided calculations (see \S \ref{SecHC}).  
Most notably they showed (see \S \ref{InvGC} below) that 
whenever $A_\theta$ is realized as the endomorphism
algebra of a Heisenberg $A_{\theta'}$-module endowed with the
$A_{\theta'}$-valued Hermitian structure obtained by twisting the
canonical one by a positive invertible element in $A_\theta$, the
curvature of $A_\theta$ with respect to the corresponding spectral
triple over $A_{\theta'}$ is equal to the modular curvature associated
to the same element of $A_\theta$ viewed as conformal factor. 
In a certain sense
this is reminiscent of Gauss's
Theorema Egregium ``If a curved surface is developed upon any
other surface whatever, the measure of curvature in each point remains
unchanged''.
 
The fundamentals of Connes' pseudodifferential calculus as well as its extension
to twisted $C^*$-dynamical systems, which provide
the essential device for proving all the above results, are explained in \S \ref{psi-do}.
Finally, \S \ref{ResExp} clarifies how to use the affiliated symbol calculus in order
to compute the resolvent trace expansion, or equivalently the heat trace expansion, for
the relevant Laplace-type operators.
 
 
\section{Curvature of modular spectral triples}
\label{MLSec2}

\subsection{Flat spectral triples} \label{flatST}
In noncommutative geometry a metric structure on a space with
$C^*$-algebra of coordinates $A$ is represented by a triad of data
$(A, \cH, D)$ called \textit{spectral triple}, modeled on the Dirac
operator on a manifold: $A$ is realized as a norm-closed subalgebra of
bounded operators on a Hilbert space $\cH$,  $D$ is an unbounded
self-adjoint operator whose resolvent belongs to any $p$-Schatten
ideal with $p > d$ where $ d > 0$ signifies the dimension, and $D$
interacts with the coordinates by having bounded commutators (or more
generally bounded twisted commutators) $[D, a]$ for any $a$ in a dense
subalgebra of $A$.  The Dirac operator was chosen as model since it
represents the fundamental class in $K$-homology and at the same time
plays the role of a quantized inverse line element
(see~~\cite{Con2013}).  In the case of $A_{\theta}$ one can obtain
such a triad by simply reproducing the construction of the 
$\bar\pl + \bar\pl^*$ operator on the ordinary torus 
$ \T^2 = \left( \R/2\pi\Z\right)^2$.

To fix the notation we briefly review some basic properties of
the $C^*$-algebra $A_{\theta}$ with $\theta \in \R \setminus \Q$. 

First of all, the torus $ \T^2 $ acts on $A_{\theta}$ via
the representation by automorphisms defined on the basis elements by 
\begin{align} \label{actT}
  \alpha_{\textup{r}} (U_1^n \, U_2^m)=e^{i(r_1n + r_2m)}U_1^n \, U_2^m, \quad \,
      \textup{ r } = (r_1, r_2) \in \R^2 .
\end{align}
By analogy with the action of $ \T^2 $ on $A_0 \equiv C(\T^2)$, we call these automorphisms 
\textit{translations}.

The smooth vectors of this representation of $ \T^2 $ are precisely
the elements of the form 
$a= \sum_{m,n\in \mathbb{Z}}a_{m,n}U_1^m U_2^n$ with rapidly decaying
coefficients, i.e. such that 
$(1+|m|)^k (1+|n|)^\ell |a_{m,n}| \leq C_{k,\ell}$, $ k, \ell > 0$.
These elements form a subalgebra $\At$ which is the analogue of $C^\infty (\T^2 )$ 
viewed in Fourier transform.
 The assignment 
\begin{align*}
\At \ni a=\sum_{m,n\in \mathbb{Z}}a_{m,n}U_1^mU_2^n \quad
 \mapsto  \quad \vp_0(a) = a_{0,0} ,
\end{align*}
determines the unique normalized trace $\vp_0$ of the $C^*$-algebra $A_{\theta}$.

The image of the differential of the above representation 
on $\At $ is the Lie algebra generated by the outer derivations 
 $ \delta_1$ and $\delta_2 $, uniquely determined by the relations
$\delta_i(U_j)= \delta_i^jU_j $, $i,j \in \{1,2\}$.

By analogy with the ordinary torus, one defines on $A_{\theta}$ a
translation invariant complex structure with modular parameter
$\tau \in \C$, $\Im\tau >0$,
by means of the pair of derivations
\begin{align} \label{delta12}
 \delta_\tau = \delta_1 + \bar{\tau} \delta_2,  \quad \delta_\tau^\ast =\delta_1 + \tau \delta_2 ;
 \end{align}
these are the counterparts of the operators 
$ \frac 1i \left(\partial / \partial x + \bar{\tau} \partial / \partial y\right)$, and
$ \frac 1i \left(\partial / \partial x + \tau \partial / \partial y\right)$ acting on $C^\infty (\T^2)$.

To obtain the analogue of the corresponding flat metric on $\T^2$, we
 let $\cH_0 \equiv L^2(A_\theta, \vp_0)$ denote the Hilbert space completion of 
$A_{\theta}$ with respect to the scalar product
\begin{equation*}
 \langle a,b \rangle = \varphi_0(a^*b),\quad a,b \in \At \, .
\end{equation*}
The space $\Omega^1 \At $ of formal  $1$-forms 
$\sum \, a \ d b $,  $a,b \in \At $,
is also endowed with a semi-definite inner product defined by  
\[
  \langle a db , a' db' \rangle = \vp_0 (a^* \,(a')
  \delta_\tau (b') \, \delta_\tau b^*) \, , \quad a,a',b,b' \in \At .
\] 
On completing its quotient modulo the subspace of those elements 
$\omega \in \Omega^1 \At$ such that $\langle \omega, \omega \rangle = 0$,
 one obtains a Hilbert space denoted $\cH^{(1,0)}$. 
$\cH^{(1,0)}$ is also an $A_\theta$-bimodule under
the natural left and right action of $A_\theta$. Both these actions are unitary.
Moreover, the linear map from $\Omega^1 \At $ to
$\At $ defined by sending the class of $ \sum a db$ in $\cH^{(1,0)}$  
to $\sum a  \delta_\tau(b) \in \cH_0 $
induces an $A_\theta$-bimodule isomorphism between $\cH^{(1,0)}$ and $\cH_0$. 

Denoting by $\partial_\tau $ the closure of the
operator $ \delta_\tau: \At \to \cH_0$ viewed as unbounded operator from 
$\cH_0$ to $\cH^{(1,0)}$ one obtains a spectral triple
$(A_\theta, \tilde{\cH}, D_\tau)$ by taking $\tilde{\cH} = \cH_0 \oplus \cH^{(1,0)}$ 
and as unbounded self-adjoint operator $ D_\tau= \left(
  \begin{array}{cc}
    0 & \partial_\tau^*\\
    \partial_\tau & 0\\
  \end{array}
\right) $.
Concurrently, the triad $(A_\theta^{\rm op}, \tilde{\cH}, D_\tau)$ is a spectral triple with respect
to the right action of $A_\theta$. One can turn it into a spectral triple for the left action of 
$A_\theta$ by passing to its transposed $(A_\theta, \bar{\cH}, \bar{D}_\tau)$
(see \cite[\S 1.2]{CM2014} for the general definition), where $\bar{\cH}$ is the complex conjugate
of $\tilde{\cH}$ and $ \bar{D}_\tau= \left(
  \begin{array}{cc}
    0 & \bar\partial_\tau^*\\
    \bar\partial_\tau & 0\\
  \end{array}
\right) $.

\subsection{Modular spectral triples} \label{ModST}
To implement the analogue of a conformal change of metric structure, we choose a 
self-adjoint element $h=h^* \in \At$ and use it to replace the trace 
$\varphi_0$ by the positive linear functional $\varphi \equiv \varphi_h$ defined by
\begin{align} \label{vph}
\vp(a) \equiv \varphi_h (a) = \vp_0 (ae^{-h}) \, , \quad a \in A_{\theta} \, .
\end{align}
Then $\vp$ determines an inner product $ \langle \ , \ \rangle_\varphi$ on $A_{\theta}$,
\[
\langle a,b\rangle_\varphi = \varphi (a^* b) \, , \quad a,b \in A_{\theta} \, ,
\]
which by completion gives rise to a Hilbert space $\cH_\varphi$.
The latter is again an $A_\theta$-bimodule but, since 
 $\vp$ is no longer tracial, the right action is no longer unitary.
 
The non-unimodularity of  $\vp$ is expressed by Tomita's modular operator $\Delta$,
which in this case is
\begin{equation*}
\modu (x) = e^{-h} x e^h  , \quad x \in A_{\theta} ,
\end{equation*}
and gives rise to the $1$-parameter group 
 of inner automorphisms
\begin{equation} \label{sigmat}
\sigma_t (x) \, = \, \modu^{-it} \, = \, e^{ith} x e^{-ith} , \quad x \in A_{\theta} .
\end{equation}
Instead of the tracial property  $\vp$ satisfies the KMS condition
 \begin{equation*}
\varphi (ab) = \varphi (b \sigma_i \, (a)) = \varphi (b e^{-h} ae^{h}) \, , \quad a, b \in A_{\theta} .
\end{equation*}
 To restore the unitarity of the right action one redefines it by setting 
\[
  a^{\rm op} := J_\varphi a^* J_\varphi \in \cL (\cH_\vp), \quad  a \in A_{\theta} ,
\]
 where $ J_\vp (a) \, = \Delta^{1/2}( a^\ast)=\, k^{-1} a^\ast k$, $a \in A_{\theta}$,
and  $k=e^{h/2}$.
 
While keeping $\cH^{(1,0)}$ unchanged, we now view $\delta_\tau$ as a densely
defined operator from $\cH_\vp$ to $\cH^{(1,0)}$. Its closure 
$\pl_\vp$ is then used to define 
$ D_\vp= \left(
  \begin{array}{cc}
    0 & \partial_\vp^*\\
    \partial_\vp & 0\\
  \end{array}
\right) $ 
giving rise to the triad $(A_\theta^{\rm op},\tilde{\cH}_\vp ,D_\vp)$, where
$\tilde{\cH}_\vp = \cH_\vp \oplus \cH^{(1,0)}$. 
This is a \textit{twisted} spectral triple
(see~\cite{cmIII} for the general definition) over $A^{\rm op}$, 
with the twisted commutators  
$ D_\vp \,a^{\rm op} - (k^{-1}ak)^{\rm op} D_\vp $, \,
$a \in \At$  bounded. Its {\em transposed}, formed as in the flat case, yields
the \textit{modular spectral triple} over $\At$, with operator
$\bar{D}_\vp =  \left(
  \begin{array}{cc}
    0 & k \pl_\vp\\
     \pl_\vp^*k & 0\\
  \end{array}
\right) $, where the conformal factor $k$ acts by left multiplication, and with
underlying Hilbert space $\overline{\tilde{\cH}}_\vp$.

By a series of identifications, it is shown in~\cite[\S 1.3]{CM2014} that  
the modular spectral triple associated to $\vp$, or equivalently to the
conformal factor  $k=e^{h/2}$, is canonically isomorphic to the twisted spectral triple
$(A_\theta , \tilde{\cH}_0,  D_k)$ with $\tilde{\cH}_0 := \cH_0 \oplus \cH_0$ and
$D_k : =  \left(
  \begin{array}{cc}
    0 & k \delta_\tau\\
     \delta_\tau^*k & 0\\
  \end{array}
\right) $. 

We finally note that $D_k^2 =  \perta \oplus \pertb$, where 
\begin{align} \label{pertab}
 \perta := k \lapa_\tau k \equiv
k \delta_\tau   \delta_\tau^* k \quad \text{and} \quad \pertb = \delta_\tau^* k^2 \delta_\tau ,
\end{align}  
are the counterparts of the Laplacian on functions, respectively the
Laplacian on $(0,1)$-forms.

\subsection{Modular curvature} \label{ModC}
The meaning of locality in noncommutative geometry is guided by the analogy
with the Fourier transform, which interrelates the local behavior of functions with
the decay at infinity of their coefficients. In a similar way, in the noncommutative 
formalism the local invariants of a spectral triple $(A, \cH, D)$ are encoded in 
the high frequency behavior of the spectrum of the 
`inverse line element' $D$ coupled with the action of the algebra of coordinates. 
For example, the local index formula in noncommutative geometry~\cite[Part II]{ConMos:LIF}
expresses the Connes-Chern character 
of a spectral triple with finite dimension spectrum in terms
of multilinear functionals given by residues of zeta functions defined by
 \begin{align*}
  z \mapsto \Tr \left(a_0 [D, a_1]^{(k_1)} \ldots [D, a_p]^{(k_1)} \, |D|^{-z} \right),   
  \quad \Re(z) >>0 ,
\end{align*}
where $ a_0, \ldots, a_p \in \sA$ and
$[D, a]^{(k)} = [D^2, \ldots , [D^2, [D,a]] \cdots ]$ with $D^2$ repeated $k$-times;
the existence of the meromorphic continuation of such zeta functions is built in the
definition of  \emph{finite dimension spectrum} for a spectral triple.
Clearly, perturbing $D$ by a trace class operator will not affect these residue
functionals, whence the local nature of the index formula described in their terms.
 
In the specific case of the noncommutative torus 
the concept of locality can be pushed much closer to the customary one. 
Namely, if $(A_\theta , \tilde{\cH}_0,  D_k)$ is a modular spectral
triple as in \S \ref{ModST},
for its Laplacian `on functions' there is an asymptotic expansion 
\begin{align} \label{heat-asy}
  \Tr \bl a\, e^{-t\perta} \br\, \sim_{t \searrow 0} \,
   \sum_{q=0}^\infty  {\rm a}_{2q}(a, \perta)\, t^{q-1} , \qquad \, a \in \At ,
\end{align} 
whose functional coefficients $ {\rm a}_{2q}$ are not only local in
the above sense, but they are also absolutely continuous with respect
to the unique trace, i.e. of the form
\begin{align*}
\At \ni a \longmapsto {\rm a}_{2q}(a, \perta) =  \vp_0 (a \, \cK^{(q)}_k) , \qquad
\cK^{(q)}_k \in \At ,
\end{align*}
with `Radon-Nikodym derivatives' $\cK^{(q)}_k \in \At $ computable 
by means of symbolic calculus. The technical apparatus which justifies
the heat expansion \Eqref{heat-asy} as well as the explicit computation of
$\cK^{(0)}_k$ will be discussed in  \S\S \ref{psi-do}-\ref{ResExp}.

In particular, the Radon-Nikodym derivative of the term $ {\rm a}_2$,
which classically delivers the scalar curvature, was fully computed
in~\cite{CM2014, FarKha2013} and represents the 
modular scalar curvature. Abbreviating its notation to $ \cK_k$
instead of $\cK^{(0)}_k$, it has the following expression:
\begin{equation}\label{a2term}
 \cK_k=-\frac{\pi}{2 \Im\tau}\left(K_0(\nabla)(\lapa(h))+
 \frac 12 H_0(\nabla^{(1)},\nabla^{(2)})(\square_\Re (h)\right) ,
\end{equation}
where $\nabla=\log \modu$ is the inner derivation implemented  by $-h$,
$$
 \lapa(h)= \delta_\tau \delta_\tau^* =
  \delta_1^2(h)+2 \Re\tau \, \delta_1\delta_2(h)+|\tau|^2 \delta_2^2(h) ,
 $$
 $\square_\Re$ is the Dirichlet quadratic form
 $$
 \square_\Re (\ell) :=
(\delta_1(\ell))^2+ \Re\tau\left(\delta_1(\ell)\delta_2(\ell)+\delta_2(\ell)\delta_1(\ell)\right)+|\tau|^2 (\delta_2(\ell))^2\,,
$$
and $\nabla^{(i)}, \, i=1, 2$, signifies that $\nabla$ is acting on the $i$th factor.
The functions $K_0(s)$ and $H_0(s,t)$, whose expressions resulted from the
symbolic computations, are given by 
  \begin{align*} 
  \begin{split}
 &  \qquad  \qquad K_0(s)=\frac{-2+s\, {\rm coth}\left(\frac{s}{2}\right) }{s\, \sinh\left(\frac{s}{2}\right)}
 \qquad \text{and}  \qquad   H_0(s,t)= \\
&\frac{t (s+t) \cosh(s)-s (s+t) \cosh(t)+(s-t) (s+t+\sinh(s)+\sinh(t)-\sinh(s+t))}{s t (s+t) \sinh\left(\frac{s}{2}\right) \sinh\left(\frac{t}{2}\right) \sinh\left(\frac{s+t}{2}\right)^2}.
\end{split}
\end{align*}
The second function is related to the first by the functional identity  
\begin{align}\label{FI}
 - \frac 12  \tilde H_0(s_1,s_2)& = \frac{\tilde K_0(s_2)-\tilde K_0(s_1)}{s_1+s_2} +\\ \notag
  &\frac{\tilde K_0(s_1+s_2)-\tilde K_0(s_2)}{s_1}-\frac{\tilde K_0(s_1+s_2)-\tilde K_0(s_1)}{s_2}  ,
\end{align}
where
\begin{align} \label{tildeKH}
  \begin{split}
\tilde{K}_0(s)\, &= \, 4\frac{\sinh(s/2)}{s} K_0(s) \qquad \text{and} \\
\tilde{H}_0(s,t)\, &= \, 4\frac{\sinh((s+t)/2)}{s+t} H_0(s,t) .
\end{split}
\end{align}
A noteworthy feature of the main curvature-defining function is that,
up to a constant factor,  $\tilde K_0$ is a generating function for the Bernoulli numbers; 
precisely,
\begin{equation} \label{Bern}
 \tilde K_0(t)= 8\sum_{n=1}^\infty \frac{B_{2n}}{(2n)!}t^{2n-2} .
\end{equation}
 
\subsection{Modular Gauss-Bonnet formula}   \label{ModGB}
Since the $K$-groups of the noncommutative torus are the same as of the ordinary
torus, its Euler characteristic vanishes. Thus, the analogue of the Gauss-Bonnet theorem
for the modular spectral triple is the identity
\begin{align*}
\vp_0 \Bl  \cK^{(0)}_k \Br \, = \, 0 .
\end{align*}
This can be directly checked by making use of the fact that the group of modular automorphisms
$\sigma_t$ (cf. Eq. \eqref{sigmat}) preserves the trace $\vp_0$ and fixes the dilaton $h$,
in conjunction with the `integration by parts' rule
\begin{align*}
\vp_0(a \delta_j(b)) \, = \, - \vp_0 (\delta_j(a) b) , \qquad a,b \in \At \, .
\end{align*}
(See~\cite[Lemma 4.2]{CM2014} for the precise identity to be used). 

An alternative variational argument, given in \cite[\S 1.4]{CM2014}), runs as follows.  
Consider the
family of Laplacians 
\begin{align} \label{Ds}
 \triangle_s := \, k^s\,\lapa\,k^s = e^{ \frac{s h}{2} }\,\lapa_\tau\,e^{\frac{s h}{2} }, \qquad s \in \R.
\end{align}  
One has $\displaystyle
 \frac{d}{ds}\triangle_s =  \frac{1}{2} \left(h \triangle_s  + \triangle_s  h\right)$.
By Duhamel's formula one can interchange the derivative with the trace. hence 
\begin{eqnarray*}
\frac{d}{d s} \Tr \bl e^{-t \triangle_s } \br  \, = \, -  t  \, \Tr \bl h \, \triangle_s \,
e^{-t \triangle_s } \br
\, = \, t  \,\frac{d}{dt}  \Tr \bl h \, e^{-t \triangle_s } \br .
 \end{eqnarray*}
Differentiating term-by-term the asymptotic expansion Eq.~\eqref{heat-asy} 
(with $a=1$ omitted in notation) yields
   \begin{equation*}  
\frac{d}{d s} {\rm a}_{j}(\triangle_s )  \, = \,
\frac{1}{2} (j-2) \, {\rm a}_{j}(h, \triangle_s ) \, , \quad j \in \mathbb{Z}^+ .
\end{equation*}
In particular, ${\rm a}_{2}(\triangle_s ) = {\rm a}_{2}(\lapa_\tau)$. The latter
vanishes because $\lapa_\tau$ is isospectral to the Laplacian of the ordinary torus 
with the same complex structure and, as is well-known, if $\lapa_M$ is the
Laplacian on a Riemann surface then ${\rm a}_{2}(\triangle_M) = \frac{\chi(M)}{6}$,
where $\chi(M)$ is the Euler characteristic of $M$.

\subsection{Variation of determinant and modular Gaussian curvature} \label{var-det}

The  zeta function
$\zeta_{\perta} (a, z) = \Tr \bl a\, \perta^{-z}\, (1-P_k) \br$ ,  $\Re z >  2$ 
where $P_k$ stands for the orthogonal projection onto $\Ker \perta$, is
related to the corresponding theta function by the Mellin transform
\begin{align*} 
 \zeta_{\Lapl_k} (a, z) \, = \, \frac{1}{\Gamma (z)} \int_0^\infty t^{z-1} \,
  \Tr \Bl a\, \bl e^{-t \Lapl_k} - P_k \br\Br  \, dt\, .
\end{align*}
The asymptotic expansion Eq.~\eqref{heat-asy} ensures that it has meromorphic continuation
and its value at $0$ is
\begin{align}   \label{z0}
 \zeta_{\Lapl_k}(a, 0) =  {\rm a}_2 (a, \triangle_\vp)  -  \Tr (P_k \, a \, P_k)  
 =  {\rm a}_2 (a, \triangle_\vp)  -  \frac{\vp_0(ak^{-2})}{\vp_0(k^{-2})} .
\end{align}  
In particular for $a=1$ (suppressed in notation), one has
\begin{align} 
 \zeta_{\Lapl_k}(0) \, = \, -1 ,
 \end{align}
and also the Ray-Singer log-determinant is well-defined:
\begin{align*}  
 \log \Det \Lapl_k := -  \zeta'_{\Lapl_k} (0) .
\end{align*}
 
Differentiating the $1$-parameter family of zeta functions corresponding to \eqref{Ds}
 one obtains the identity
\begin{align*} 
\frac{d}{ds} \zeta_{\triangle_{sh}} (z)\, = \,
- z \, \zeta_{\triangle_{sh}} (h, z) , \quad \forall \, z \in \mathbb{C} ,
\end{align*}
which in turn yields the variation formula
\begin{align*} 
- \frac{d}{ds} \zeta^\prime_{\triangle_{sh}} (0)\, =
 \, \zeta_{\triangle_{sh}} (h, 0)  .
\end{align*}
From Eq.~\eqref{z0} and Eq.~\eqref{a2term} applied to the weights $\vp_s$ with dilaton $sh$
one obtains
\begin{align*}  
\log \Det \perta&= \,\log \Det \triangle  +  \log \vp (1)  
 - \frac{\pi}{\Im\tau} \int_0^1 \vp_0\bigl( h\big(s K_0(s\nabla)(\triangle (\log k)) \\
&+s^2 H_0(s\nabla^{(1)},s\nabla^{(2)})(\square_\Re(\log k))\big)\bigr) ds
\end{align*}
The first term is the same as for the corresponding elliptic curve and
by the Kronecker limit formula has the expression (cf. \cite[Theorem~4.1]{RaySin73})
\begin{align*}
 \log \Det \triangle  \, =- \frac{d}{ds} \big|_{s=0} \sum_{(n,m)\neq (0,0)}|n+m\tau|^{-2s}=\, 
\log \left(4 \pi^2 \,  |\eta(\tau)|^4 \right),
\end{align*}
where  $\eta$ is the Dedekind eta function $\eta (\tau) = e^{\frac{\pi \, i}{12} \, \tau} 
 \prod_{n>0} \left(1-e^{2\pi i n \tau} \right)$.
After a series of technical manipulations of the last term (see~\cite[\S 4.1]{CM2014}), 
one obtains the {\em modular analogue of Polyakov's anomaly formula}:
\begin{equation} \label{-F}
 \log \Det \perta = \log \left(4 \pi^2 \,  |\eta(\tau)|^4 \right) +  \log \vp (1)
   - \frac{\pi}{4\Im\tau}  \vp_0\left(K_+ (\nabla^{(1)} )(\square_\Re(h))\right) ,
 \end{equation}
where $K_+(s) := \frac{4}{s^2}-\frac{2 {\rm coth}(\frac{s}{2})}{v} \geq 0$, $s \in \R$. 
Furthermore, it is
shown in~\cite[Proof of Theorem 4.6]{CM2014} that the positivity of the function $K_+$
can be upgraded to operator positivity, implying the inequality
\begin{equation} \label{pos}
 \vp_0\left(K_+ (\nabla^{(1)} )(\square_\Re(\log k))\right) \geq 0 ,
\end{equation}
with equality only for $k=1$. 

The (negative of)
log-determinant can be turned into a scale invariant functional by adding the area term:
\begin{align} \label{defF}
F (\log k): = \, \zeta'_{\Lapl_k}(0)+ \log \vp (1) \ = \, - \log \Det (\Lapl_k) + \log \vp (1) .
\end{align}
Due to the equality Eq.~\eqref{z0}, the corrected functional $F$ remains unchanged
when the Weyl factor $k$ is multiplied by a scalar.
In the new notation the identity Eq.~\eqref{-F} reads as follows:
\begin{equation} \label{+F}
F(h) = - \log \left(4 \pi^2 \,  |\eta(\tau)|^4 \right)  
   + \frac{\pi}{4\Im\tau}  \vp_0\left(K_+ (\nabla^{(1)})(\square_\Re(h))\right) .
 \end{equation}
In view of the inequality Eq.~\eqref{pos} one concludes
that, as in the case of the ordinary torus (cf.~\cite{OPS88}),
 {\em the scale invariant functional $F$
attains its extremal value only for the trivial Weyl factor}, in other words
at the flat metric.

The gradient of $F$ is defined by means of the inner product of $L^2(A_\theta, \vp_0)$
via the pairing
\begin{align*}
\langle \grad_h F , a \rangle  \equiv \, \vp_0 (a\,  \grad_h F )  
= \, \frac{d}{d\varepsilon}\big|_{\varepsilon=0}
F (h + \varepsilon a) , \quad
  a = a^\ast \in \At .
\end{align*}
A direct computation of the gradient, using the definition Eq.~\eqref{defF}
combined with the identities Eq.~\eqref{z0} and Eq.~\eqref{a2term}, yields the following
explicit expression (cf.~\cite[Theorem 4.8]{CM2014}):
 \begin{align} \label{grad}
   \grad_h F = { \frac{\pi}{4\Im(\tau)}\left( \tilde K(\nabla)(\Lapl(h))+
 \tilde H(\nabla^{(1)},\nabla^{(2)})(\square_\Re (h))\right) } .
 \end{align}
 In the case of the ordinary torus the gradient of the corresponding functional 
 (cf.~\cite[(3.8)]{OPS88}) gives precisely the Gaussian curvature. This 
makes it compelling to take the above formula as definition of the
{\em modular Gaussian curvature}.
 
Finally, computing the gradient of $F$ out of its explicit formula Eq.~\eqref{+F}, and 
then comparing with the expression Eq.~\eqref{grad}, produces the 
functional identity Eq.~\eqref{FI} relating $\tilde H$ and $\tilde K$.

\section{Morita invariance of the modular curvature} \label{MorInv}

\subsection{Foliation algebras and Heisenberg bimodules} \label{FolAlg}
The most suggestive depiction of the noncommutative torus was given by Connes
in~\cite{Con1982}, where he described it as the ``space of leaves'' 
for the Kronecker foliation $\cF_\theta$ 
of the ordinary torus $\T^2 = \R^2/\Z^2$, given by the 
differential equation $dy -\theta dx =0$ with $\theta \in \R \setminus \Q$. 
The holonomy groupoid  $\cG_\theta$ 
of this foliation identifies with the smooth groupoid determined by
the flow of the above equation. Its convolution $C^*$-algebra 
 $C^*(\cG_\theta)$, which represents the (coordinates of the) space of leaves, coincides
 with the crossed product $C(\T^2) \times_\theta \R$, where the action of $\R$ on $\T^2$  
 is given by the flow \Eqref{actT}. 
 $C^*(\cG_\theta)$ is isomorphic to $ A_\theta \otimes \cK$, where $\cK$ denotes 
 the  $C^*$-algebra of compact operators, and thus strongly Morita equivalent to  $A_\theta$. 
   
Finer geometric representations of the space of leaves are obtained by passing to reduced 
$C^*$-algebras associated to complete transversals. Any pair of relatively prime integers 
$(d, c) \in \Z^2$ determines a family of lines of slope $\frac{d}{c}$, which project onto simple
closed geodesics in the same free homotopy class, and the free homotopy classes of closed
geodesics on $\T^2 $ are parametrized by the rational projective line 
$P^1(\Q) \equiv \Q \cup\{\frac{1}{0}\}$. Letting $N_{c, d}$ denote the
primitive closed geodesic of slope $\frac{d}{c}$ passing through the base
point of $\T^2$, one obtains a complete transversal for $\cF_\theta$. The convolution
algebra of the corresponding  \'etale holonomy groupoid identifies
with the crossed product algebra $ C(\R/\Z) \rtimes_{\gt'}\Z$, where $1 \in \Z$
acts by the rotation of angle  $\gt' =\frac{a\gt +b}{c \gt +d}$ with
$a, b \in \Z$ chosen such that $ad-bc=1$. This $C^*$-algebra is none other than $A_{\gt'}$.
In particular $A_\gt = C(\R/\Z) \rtimes_{\gt}\Z$ is the reduced $C^*$-algebra associated to
$N_{0, 1}$. 
By construction all algebras $A_{\gt'}$ with $\gt' =g \cdot \gt$, $g \in \SL(2,\Z)$, are
Morita equivalent, and they actually exhaust (cf. \cite{Rie1981}) all the
noncommutative tori Morita equivalent to $A_\gt$. 

In the same framework Connes~\cite[\S 13]{Con1982} gave a geometric description of the   
$(A_{\theta'}, A_\theta)$-bimodules $E(g,\theta)$ implementing the Morita 
equivalence of $A_\gt$ with $A_{\gt'}$.  $E(g,\gt)$ is a completion
of the  $(\Atp, \At)$-bimodule
$ \sE(g,\gt):= \sS(\R)^{|c|} \equiv \sS(\R\times \Z_c)$, $\Z_c:=\Z/c\Z$, with
the actions defined as follows:
\begin{align*}
\begin{split} 
&(fU_1)(t,\ga):= e^{\tpii (t- \frac{\ga d}{c})} f(t,\ga)\, , \,
        (fU_2)(t,\ga):= f(t-\frac{c\gt + d}{c} ,\ga-1) ; \\
& (V_1f)(t,\ga):= 
          e^{\tpii \bigl(\frac{t}{c\gt+d}- \frac{\ga}{c}\bigr)} f(t,\ga)\, , \,
        (V_2f)(t,\ga):= f(t-\frac 1c,\ga-a).
\end{split}
\end{align*}
If $c=0$ then  $E(g,\gt)= A^{\rm op}_\theta $ is the trivial 
$(A_{\theta}^{\rm op}, A_\theta)$-bimodule.
By analogy with the vector bundles over elliptic curves, one defines the rank, degree 
and slope of  $\sE(g,\gt)$ by
${\rk} (g,\theta) = c\theta +d$,  $\deg(g,\theta) = c$, resp.
$\mu (g,\theta):=\frac{\deg(g,\theta)}{\rk(g,\theta)}$. 
 
The $L^2$-scalar product on $E(g,\gt)$
\begin{align*}
  <f_1,f_2>:= \int_{\R\times \Z_c} \ovl{f_1(t,\alpha)} f_2(t,\alpha) dt d\alpha 
\end{align*}
where the integration is with respect to the
Lebesgue measure on $\R$ and the counting measure on $\Z_c$, determines
uniquely $A_\theta$--valued and $A_{\theta'}$--valued  inner products satisfying
the double equality
\begin{align} \label{Herminn}
|{\rk} (g,\theta)|\, \varphi'_0( _{A_{\theta'}}<f_2,f_1>)\, = \, <f_1,f_2>\, =\, 
 \varphi_0(<f_1,f_2>_{A_\theta}) ,
\end{align}
where $ \varphi'_0$ stands for the trace of $A_{\theta'}$.
 The completion $E(g,\theta)$ with respect to 
$|<\cdot,\cdot>_{A_\theta}|^{1/2}$ is a full right $C^*$--module over $A_\theta$,
and $\End_{A_\theta} = A_{\theta'}$. In addition,
 $ \cH_0(g,\theta) := E(g,\theta)\otimes_{A_\theta} L^2(\mathcal{A}_\theta,\varphi_0)$  
is the Hilbert space $L^2(\R\times \Z_c)$.

Instead of the $\R^2$--action Eq.~\eqref{actT}, the non-trivial bimodules $\sE(g,\gt)$
are acted upon by the Heisenberg group $H_3(\R)$.
Equivalently,  $\R^2$ acts projectively on $\sE(g,\gt)$, and this action is
compatible with the natural $\R^2$--actions on $\At$ and $\Atp$.
At the Lie algebra level, this action gives rise to the standard connection  
$\nabla^\sE$ on 
$\sE(g,\gt)$, given by the derivatives 
$(\nabla_1 f)(t,\ga) = \frac{\pl}{\pl t} f(t,\ga),$
$ (\nabla_2 f)(t,\ga) = {\tpii}\, t\, \mu (g,\gt) \, f(t,\ga)$ with
constant curvature: $[\nabla_1,\nabla_2]=\,{\tpii} \mu (g,\gt)\, \Id$. Furthermore,
this connection is bi--Hermitian, in the sense that it preserves both the
$A_\theta$--valued and the $A_{\theta'}$--valued inner product.

 \subsection{Modular Heisenberg spectral triples} \label{ModHST}
Each bimodule $\sE=\sE(g,\theta)$ gives rise to a double spectral
triple, by coupling it with the flat Dirac $D_\tau$ by means of its
standard connection. Specifically, $\nabla^\sE$ splits into holomorphic
and anti-holomorphic components, $\nabla^\sE= \pl_\sE \oplus \pl_\sE^*$,
 where  $\pl_\sE:=\nabla_1+\taubar\nabla_2 $. One then forms
the operator  $D_\sE= \left(
  \begin{array}{cc}
    0 & \partial_\sE^*\\
    \partial_\sE & 0\\
  \end{array} \right)$ acting on the Hilbert space $\tilde{\cH}(g,\gt)
 = \cH_0(g,\gt) \oplus \cH^{(1,0)}(g,\gt)$, where
 $\cH^{(1,0)} (g,\gt) := E(g,\gt)\otimes_{\At}\cH^{(1,0)}(\At)$.
Together with the natural right action of $A_\gt$, these data define
a spectral triple of constant curvature
$\left(A_\gt^{\rm op}, \tilde{\cH}(g,\gt), D_\sE\right)$. We note that 
from the spectral point of view
the operator $D_\sE$ resembles the Hodge-de Rham operator
of an elliptic curve with coefficients in a line bundle. In particular its  
Laplacian $ \Lapl_\sE = \pl_\sE^*\pl_\sE $
is a direct sum of $|\deg (\sE)|$ copies of the harmonic oscillator
\begin{align*} 
H :=  
 - \frac{d^2}{d t^2} + 4 \pi^2 \mu(\sE)^2 |\tau|^2 t^2 
 - 4 \pi i \mu(\sE) \Re(\tau) \, t\, \frac{d}{d t}- \tpii \mu(\sE) \taubar \Id . 
\end{align*}

Now turning on the conformal change Eq.~\eqref{vph} from $\vp_0$ to $\vp_h$,
one replaces  $D_\sE $ by $D_{\sE, \vp}$ in the same way as in \S \ref{ModST}. 
The resulting
spectral triple over the algebra $A_\gt^{\rm op}$ is again a twisted one.
After correcting for the lack of unitarity of the action of $A_\gt^{\rm op}$ again 
as in \S \ref{ModST}, the operator $D_{\sE, \varphi}$ is being canonically identified
with $D_{\sE, k}:= \begin{pmatrix}
  0   &   R_k\partial_{\sE}^*  \\ 
\partial_{\sE}R_k & 0
\end{pmatrix}$  acting on 
$\tilde{\cH}_0(g,\gt) = \cH_0(g,\gt) \oplus \cH_0(g,\gt)$. 

The appropriate {\em transposed} in this setting is constructed using the 
canonical anti-isomorphism from $\sE =\sE(g, \gt)$ to $\sE' :=\sE(g^{-1}, \gt')$, 
\[
J_{g,\gt}(f)(x,\ga)=   \ovl{f((c\gt+d)x,\, -d^{-1}\ga)} , \quad f \in \sE(g,\gt) ,
\] 
which  switches the $(\Atp, \At)$--action on the first into the $(\At, \Atp)$--action on
the second.
We thus arrive at the modular Heisenberg spectral triple
$\left(A_\gt, \, \tilde{\cH}_0(g^{-1},\gt'), \,\ovl{D}_{\sE', k}\right)$
with  
$\ovl{D}_{\sE', k}= - \rk (\sE')
\begin{pmatrix}
  0   &   k {\partial}_{\sE'} \\
  {\partial}_{\sE'}^* k & 0 
\end{pmatrix}$.
Its Laplacian on sections is 
${\triangle_{\sE', k}\, = \rk (\sE')^2\, k \partial_{\sE'} \partial_{\sE'}^* k }$.

A moment of reflection shows that replacing $\vp_0$ by $\vp$ is equivalent to
 changing the Hermitian structure on $\sE'$ by 
$k^{-2 }\in \At \equiv \End_{\Atp}(\sE')$. Indeed, in view of
 Eq.~\eqref{Herminn} applied to $\sE'$, the passage to the $\Atp$--valued 
Hermitian inner product
\begin{align} \label{Hermink}
 (f'_1,f'_2)_{\Atp, k} = |\rk(\sE')|^{-1} <k^{-2} f'_1,f'_2>_{\Atp} , \quad
f'_1,f'_2 \in \sE' ,
\end{align}
has the same effect on the $L^2$-inner product, since
\begin{align*}
\begin{split}
\vp'_0 \left((f'_1,f'_2)_{\Atp, k}\right)&= |\rk(\sE')|^{-1}\vp'_0(<k^{-2} f'_1,f'_2>_{\Atp} ) 
 = \vp_0( _{\At} <k^{-2} f'_1,f'_2> ) \\
&= \vp_0( _{\At} < f'_1,f'_2> k^{-2} ) 
=  \vp ( _{\At} < f'_1,f'_2> ) .
 \end{split}
 \end{align*}   
In conclusion, the passage from the `constant curvature metric' on $A_\theta$
represented by the Heisenberg spectral triple 
$\left(A_\gt^{\rm op}, \tilde{\cH}(g,\gt), D_\sE\right)$
 to the `curved metric' represented by the modular Heisenberg spectral triple 
$\left(A_\gt, \, \tilde{\cH}_0(g^{-1},\gt'), \,\ovl{D}_{\sE', k}\right)$ can be interpreted
as being {\em effected by changing the Hermitian structure of $\sE'$
 according to} Eq. \eqref{Hermink}. Note that this interpretation remains valid even
when $c=0$, i.e. for $\sE = \At$.

The extended version of Connes' pseudodifferential calculus (see \S \ref{psi-do}) 
allows to establish the heat asymptotic expansion
  \begin{align} \label{heat-ext}
   \Tr \left(a\, e^{-t \triangle_{\sE', k}} \right)\, \sim_{t \searrow 0} \,
  \sum_{q=0}^\infty {\rm a}_{2q}(a, \triangle_{\sE', k})\, t^{q-1}  , \qquad a \in \At ,
\end{align}
and express its functional coefficients in local form. In particular, the
curvature functional is of the form
\begin{align*}
 {\rm a}_{2}(a, \triangle_{\sE', k} =  \frac{1}{4 \pi \Im\tau}\vp_{\sE'} (a \, \cK_{{\sE', k}})
 = \frac{1}{4 \pi \Im\tau}{\rm rk}(\sE') \vp_0 (a \, \cK_{{\sE', k}}) , \quad a \in \mathcal{A}_\theta ,
\end{align*} 
 where $ \vp_{\sE'} \, := \, {\rm rk}(\sE') \vp_0$ is the natural trace on
 $\, \mathcal{A}_\theta =\End_{\mathcal{A}_{\theta'}}(\sE') $, and
 the curvature density has the expression (cf. \cite[Theorem 2.12]{LM2015}) 
\begin{align} \label{curv-E}
 \cK_{\sE', k} = \,
   K(\nabla)(\triangle(h))
 +H(\nabla^1,\nabla^2)\left(\square^\Re (h)\right) 
 +  \mu(\sE') 1 .
\end{align}

 \subsection{Ray-Singer determinant vs. Yang-Mills functional} \label{RSYM}

To obtain the variation formula of the Ray-Singer log-determinant functional  
\begin{align*} 
\At^{\rm sa} \ni h^*=h \longmapsto   \log \Det(\Lapl_{\cE', k}) :\, =\, 
 - \zeta^\prime_{\Lapl_{\cE', k}} (0) ,
\end{align*}
 one proceeds as in \S \ref{var-det}, starting with the insertion of the curvature
 expression \eqref{curv-E} in the derivative
$  - \frac{d}{ds} \zeta^\prime_{\Lapl_{\cE', k^s}} (0) = 
  \zeta_{\Lapl_{\cE', k^s}} (h, 0) $.
After integrating the resulting expression
one arrives at the following exact formula for the Ray-Singer determinant  
(cf. \cite[Theorem 2.15]{LM2015}) 
 \begin{align*} 
\begin{split}
 \log\Det(\Lapl_{\cE', k})&=\, \frac12 |\deg (\cE')|
   \log \big(2|\mu(\cE')| \Im(\bar\tau)\big)
   \,  - \, \frac 12 |\deg (\cE')| \, \vp_0 (h)\\
&\,- \frac{ |\rk(\cE')|}{16 \pi \Im\tau}\biggl(\frac{1}{3} \vp_0(h\Lapl h) +  \vp_0
 \Bl K_2(\nabla_h^1 )\bl \square^\Re(h)\br\Br\biggr).
\end{split}
\end{align*}
The scale invariant form of the functional is
 \begin{align*} 
F_{\cE'} (h) \, =\, - \log \Det(\Lapl_{\cE', k}) 
\, - \, \frac12 |\deg(\cE')| \vp_0 (h) .
\end{align*}
Using the preceding formula, it's exact expression is seen to be 
\begin{align}  \label{OPS-ex}
\begin{split}
F_{\cE'} (h) =&\, - \,\log \Det(\Lapl_{\cE'}) \\
&\, + \frac{ |\rk(\cE')|}{16 \pi \Im\tau}\left(\frac{1}{3} \vp_0(h\Lapl h) +  \vp_0
 \left(K_2(\nabla_h^1 )(\square^\Re(h))\right)\right) .
\end{split}
\end{align}
When viewed as a functional on the (positive cone of) metrics on the Heisenberg left
$\At$--module $\cE'$, {\em $F_{\cE'} $ attains its minimum only at the metric whose
corresponding connection compatible with the holomorphic
structure has constant curvature} (cf. \cite[Theorem 2.16]{LM2015}). 

Thus the Ray-Singer functional behaves in the same manner as the 
Yang-Mills functional of Connes and Rieffel (cf.~\cite{ConRie1987}), which however
is defined on the space of connections on the noncommutative torus.

\subsection{Invariance of the Gaussian curvature}   \label{InvGC}
 The gradient of the functional $F_{\cE'}$, now defined via the equation 
\[ 
\langle \grad_h F , a \rangle_ {\cE'} \equiv 
 \frac{1}{4\pi \Im\tau}\,\vp_{\sE'}(a\cdot  \grad_h F_{\cE'}):=
\frac{d}{d\epsilon}\big|_{\epsilon=0}F(h + \epsilon a) , 
\]
Its explicit expression can be computed as in \cite[\S 4.2]{CM2014}
and the answer turns out to be exactly the same as in the case of trivial coefficients,
cf. Eq. \eqref{grad}:
\begin{align*} 
 \grad_h F_{\cE'}=  { \frac{\pi}{4\Im(\tau)}\left( \tilde K(\nabla)(\Lapl(h))+
 \tilde H(\nabla^{(1)},\nabla^{(2)})(\square_\Re (h))\right)} =   \grad_h F.
  \end{align*}
  
 This result can be interpreted as expressing the invariance of the
 modular Gaussian curvature under Morita equivalence in two different ways.
 First it shows that {\em the Gaussian curvature associated to a 
 change of Hermitian metric on a Heisenberg equivalence bimodule $\sE'$ 
 by a fixed positive invertible $k \in \At$, viewed as an element of $ \End_{\cA_{\gt'}}(\cE')$,
  is independent of $\sE'$}. Secondly, regarding the Heisenberg spectral triples 
  with inverse line-element $D_{\sE', k}$ as {\em right spectral triples, conferring metrics
 to $\Atp$}, it proves that {\em the
 entire collection of Morita equivalent algebras $\{\sA_{g\cdot \gt} ; \, \gt \in \SL(2, \Z)\}$ 
 inherits the same modular curvature as the intrinsic one of $\At$}.

\section{Pseudodifferential multipliers and symbol calculus} \label{psi-do}

\newcommand{\PsiDO}{$\Psi$DO}
 
The main technical device which was used for proving the above results is a
pseudodifferential calculus adapted to twisted $C^*$-dynamical
systems, extending the well-known calculi due to Connes~\cite{Con:CAG}.

Originally, pseudodifferential operators (\PsiDO) where invented
(see Kohn and Nirenberg \cite{KN65} or for a textbook Shubin 
\cite{Shu:POS}) to study elliptic partial differential operators.
\PsiDO\ form an algebra which contains differential operators
\emph{and} the parametrices to elliptic differential operators. 
They come with a symbolic calculus: while the (complete) symbols
of differential operators are \emph{polynomials} in the covariables,
\PsiDO\ are obtained by allowing more general types of symbol
functions, e.g. of H\"ormander type.

\newcommand{\sSS}{\sS(\R^n,\sS(\R^n))}
\subsection{Ordinary \PsiDO\ in $\R^n$ from the point of view
of $C^*$--dynamical systems}

\subsubsection{Standard representation on the $L^2$-space (GNS space)}
Connes' pseudodifferential calculus on a $C^*$-dynamical
system $(A,\R^n,\ga)$ should be viewed as a pseudodifferential
calculus on $\R^n$. To motivate the defining formulas and to 
connect to the standard pseudodifferential calculus, 
we briefly recast the latter in the language of $C^*$--dynamical
systems such that the link becomes apparent. Recall that
for a suitably nice symbol function $\sigma(\xi,s), s,\xi\in \R^n$
\footnote{For consistency with the later exposition
we deliberately use a somewhat unusual order and naming
convention for the variables $\xi,s$. $\xi$ plays the role
of the covariable and the spacial variable $s$ is normally
called $x$ in \PsiDO\ textbooks.

We want to view the function $\sigma(\xi,\cdot)$ as an
  algebra valued function on $\R^n_\xi$. }
one defines the \emph{pseudodifferential operator with complete symbol
$\sigma$} as
\begin{equation} \label{EqML1}
  \bl\Op(\sigma) u \br (s) := \int_{\R^n} e^{i \inn{s,\xi} }
     \sigma(\xi,s) \hat u(\xi) \dsl\xi
     = \int_{\R^n}\int_{\R^n} e^{i \inn{s-y,\xi} }
     \sigma(\xi,s) u(y) dy \dsl\xi. 
\end{equation}
Now let us abuse this formula a little. 
Let $\Ainf:=\sS(\R^n) \subset C_0(\R^n)$ be the Schwartz space 
viewed as a $*$--subalgebra of $C_0(\R^n)$. 
It acts on itself by left multiplication. Furthermore,
there is a one parameter group of $*$--automorphisms 
$\ga_x(f) := f(\cdot-x)$ and a one parameter family of operators
$\pi_x(f):=f(\cdot - x)$ satisfying $\pi_y a \pi_{-y} = \ga_{-y}(a),
a\in\Ainf$. This gives rise to a covariant representation of 
the dynamical system $(\Ainf,\R^n,\ga)$ on the Hilbert space
$L^2(\R^n)$ which is the GNS space of the
$\ga$--invariant tracial weight $\varphi_0(f)=\int_{\R^n} f$,
i.e. the completion of $\Ainf$ with respect to the inner product
$\inn{f,g}_{\varphi_0} = \varphi_0( f^* g ) = \int_{\R^n} \ovl{f} g$.

Now given $u\in\Ainf$ and a symbol
$\sigma\in\sS(\R^n,\Ainf)=\sS(\R^n,\sS(\R^n)) =
\sS(\R^n\times\R^n)$ we continue from \Eqref{EqML1} and compute
\begin{equation}
  \begin{split}
    \bl\Op(\sigma) u\br(s)   &=  
        \int_{\R^n} \Bl \int_{\R^n} e^{i \inn{y,\xi} } \sigma(\xi,s)
                \dsl \xi\Br u(s-y) dy \\
     & = \int_{\R^n} \bl\sF\ii_{\xi\to y} \sigma(y)\br (s) u(s-y) dy \\
     & = \int_{\R^n}  \sigma^{\vee}_{\xi\to y}(y) \pi_y u dy (s),
  \end{split}\label{EqML2}
\end{equation}
with $\sigma^{\vee}_{\xi\to y} :=\sF\ii_{\xi\to y} \sigma$.

Thus symbols in the Schwartz space $\sS(\R^n,\Ainf)$ act, after
a Fourier transform in the first variable, covariantly with respect
to the natural representation of the covariance algebra
$\sS(\R^n,\Ainf)\rtimes_\ga \R$ on the GNS space of the weight
$\varphi_0$.

We note furthermore, that for $\sigma\in\sS(\R^n,\Ainf)=\sS(\R^n\times\R^n)$ the
operator $\Op(\sigma)$ is trace class and from the calculation
\Eqref{EqML2} we see that the Schwartz kernel of $\Op(\sigma)$
on the diagonal is given by $\int_{\R^n} \sigma(\xi,\cdot) \dsl\xi$,
hence
\begin{equation}
    \Tr\bl \Op(\sigma)\br = 
        \int_{\R^n} \int_{\R^n} \sigma(\xi,s) ds\dsl \xi 
      = \varphi_0 \bl( \sigma^{\vee}_{\xi\to y}(0) \br
       = \int_{\R^n}  \varphi_0 \bl \sigma(\xi,\cdot) \dsl\xi.
  \label{EqML3}
\end{equation}

We now take the Schwartz functions $\sigma^{\vee}_{\xi\to y}$ as
basic objects. Identifying $f\in\sSS$ with $\pi(f)=\int_{\R^n} f(x) \pi_x dx$
the space $\sSS$ becomes a $*$--algebra with $*$--representation
$f\mapsto \pi(f)$ on $L^2(\R^n)$. Explicitly, 
$\pi(f)\circ \pi(g) = \pi(f * g)$ and $\pi(f^*) =\pi(f^*)$, where
\begin{equation}\label{EqConvolution}  
   f^*(x) = \ga_x\bigl( f(-x)^*\bigr), \quad
   (f*g)(x)  = \int_{\R^n} f(y) \ga_y\bigl( g(x-y)\bigr) dy,
\end{equation}
resp. with the second variables spelled out, $\sS(\R^n\times\R^n)$
becomes a $*$--algebra with involution and product given by
\begin{equation}
     f^*(x,s) = \ovl{f (-x,s-x) }\,,
         (f*g)(x,s) = \int_{\R^n} f(y,s) g(x-y, s-y) dy.
\end{equation}    

\subsubsection{Pseudodifferential multipliers}
We now lift the previous $*$--representation to a ``universal''
multiplier representation as follows:

$\sSS$ is a pre-$C^*$--module with inner product
$\inn{f,g}=\int_{\R^n} f(x)^* g(x)dx$. Put
\begin{equation} 
  (a \, f)(x) = \ga_{-x} (a) f(x)\,,  (\uU_y f)(x) =  f(x-y)\,,
  \qquad a \in \sS(\R^n).
\end{equation}
Since $\uU_x a \uU_{-x} = \ga_x(a)$
this gives rise to a covariant representation of the $*$--algebra
$\sSS$ by associating to $f\in\sSS$ the multiplier 
$M_f=\int_{\R^n} f(x) \uU_x dx$.

If $\varphi$ is a $\alpha$--invariant trace on $\sS(\R^n)$
then the \emph{dual trace} $\psihat$ on $\sSS$ is given by
\begin{equation}\label{EqDualTrace}  
 \psihat(f)= \psi\bl f(0)\br=\int_{\R^n} \psi\bl\hat{f}(\xi)\br
 \dxi.        
\end{equation}
Note that $\dxi$ is the Plancherel measure of the dual group $(\R^n)^\wedge$
w.r.t. the duality pairing $(x,\xi)\mapsto e^{i\inn{x,\xi}}$.

In case of the trace $\varphi_0 = \int_{\R^n}$ from the previous
section the dual trace equals the trace \Eqref{EqML3} on the Hilbert
space representation $L^2(\R^n)$. This equality should be viewed as
a coincidence. In general the dual trace does not coincide with the
Hilbert space trace on a representation, resp. this depends on the
representation, see \S \ref{sec-Epo}.

By associating to $f\in\SRnA$ the multiplier $M_f=\int_{\R^n} f(x) \uU_x dx$
the space $\SRnA$ becomes a $*$--algebra. Putting $P_f := M_{{f^\vee}} $ and
allowing $f$ to be a symbol of H\"ormander class $\SmA$ we obtain an
algebra of multipliers which, via the representation $\pi$ from above,
is isomorphic to an algebra of pseudodifferential operators in $\R^n$.
We deliberately say ``an'' and not ``the'' here as in $\R^n$ there are
various versions of such algebras which differ only by the behavior
of symbols as the spacial variable $s\to\infty$, 
cf. \cite[Chap.~IV]{Shu:POS}.

\subsection{Pseudodifferential multipliers on twisted crossed products}
The action of the Heisenberg group on $\sE(g,\gt)$ induces a
$C^*$--dynamical system  $(\cA,\R^{n=2},\ga)$ ($\cA=\cA_\gt$ or $\cA =
\cA_{\gt'}$). Equivalently, $\R^n$ acts by a \emph{projective}
representation with cocycle $e(x,y):= e^{i \inn{Bx,y}},$ with a \emph{skew}
symmetric matrix $B=(b_{kl})_{k,l=1}^n$. In order to construct the
resolvent of elliptic differential operators (i.e. Laplacians)
on Heisenberg modules one therefore extends the previous
considerations to twisted $C^*$--dynamical systems. In the previous
two sections we have formulated the standard pseudodifferential
operator conventions in such a way that they carry over almost
ad verbatim to the twisted case. 

Consider a $C^*$--dynamical system $(A,\R^n,\ga)$ with, now for
simplicity, unital $A$. Furthermore, let
\begin{equation}\label{Eq1206041}  
      e(x,y):= e^{i\gs(x,y)}= e^{i \inn{Bx,y}}, \quad \gs(x,y):= \inn{Bx,y}
\end{equation}
with a skew-symmetric real $n\times n$--matrix $B=(b_{kl})_{k,l=1}^n$.
By $\sA^\infty$ we denote the smooth subalgebra, i.e. those $a\in\sA$
for which $t\mapsto \ga_t(a)$ is smooth.

As before the Schwartz space $\SRnA$ is a pre-$C^*$--module with inner product
$\inn{f,g}=\int_{\R^n} f(x)^* g(x)dx$. Putting $(\uU_y f)(x) =  e(x,-y) f(x-y)$ we obtain 
a projective family of unitaries
$\uU_x^* =\uU_{-x}$, $\uU_x\, \uU_y= e(x,y) \uU_{x+y}$,
$ x,y\in\R^n$, $\uU_x a \uU_{-x}= \ga_x(a)$, $a\in\sA^\infty$.
Together with $(a \, f)(x) = \ga_{-x} (a) f(x), a\in \sA^\infty$
and associating to $f\in\SRnA$ the multiplier $M_f=\int_{\R^n} f(x) \uU_x dx$
the space $\SRnA$ becomes a $*$--algebra.
Explicitly, cf.~\Eqref{EqConvolution}
\begin{equation}\label{EqTwistedConvolution}  
   f^*(x)    = \ga_x\bigl( f(-x)^*\bigr), \quad
      (f*g)(x)  = \int_{\R^n} f(y) \ga_y\bigl( g(x-y)\bigr) e(y,x) dy. 
\end{equation}
Note that the formula for $f^*$ is the same as in the untwisted case.

As in the untwisted case, a $\ga$--invariant trace $\psi$ on $\sA$ 
induces a dual trace $\psihat$ on $\SRnA$ which is given by the same
formula as \Eqref{EqDualTrace}.

To define the pseudodifferential operator convention we now read
\Eqref{EqML2} backwards. Namely, given Schwartz functions $f,u\in\SRnA$,
and abbreviating $f^\vee:=\sF\ii_{\xi\to y}f$ the inverse Fourier transform of $f$, 
we find
\begin{align}
  (M_{{f^\vee}}&u)(x)
         := \Bigl( \int_{\R^n} {f^\vee}(y) \uU_y u \,
                  dy\Bigr)(x)\label{EqPsiDO}\\
         &= \int_{\R^n} \ga_{-x}({f^\vee}(y)) u(x-y) e(x,-y) \, dy
         = \int_{\R^n} \ga_{-x}({f^\vee}(x-y)) u(y) e(x,y) \, dy
                  \nonumber \\
         &= \int_{\R^n}\int_{\R^n} e^{i \inn{x-y,\xi -Bx}} \ga_{-x}( f(\xi))
  u(y) dy\,\dsl\xi\label{EqPsiDOa}\\
         &= \int_{\R^n} e^{i\inn{x,\xi}} \ga_{-x}(f(\xi)) \hat u (\xi-Bx) \,
         \dsl\xi    
         = \int_{\R^n} e^{i\inn{x,\xi}} \ga_{-x}(f(\xi+Bx)) \hat u(\xi) \,
         \dsl\xi \label{EqPsiDOb}\\
         & =: (P_f u)(x)\label{eq.rev2}
\end{align}
and call the so defined multiplier $P_f$ 
a \emph{(twisted) pseudodifferential multiplier with
symbol} $f$. This should be compared to \Eqref{EqML1}.

Strictly speaking, so far we have only dealt with smoothing operators
as all symbols were Schwartz functions. One now has to extend
$P_f$ to a larger class of functions $f$. The purpose of the
somewhat lengthy exposition so far was to show that, at least in
$\R^n$ but there in a rather broad sense, smoothing operators
are nothing but convolution operators and their symbols are
obtained by applying a partial Fourier transform. General
\PsiDO\ are therefore nothing but \emph{singular} convolution
operators. This is not surprising as \PsiDO\ may,
via the Schwartz Kernel Theorem, also be viewed as singular integral
operators.

The extension to general symbol functions now follows the standard
route. Putting $P_f := M_{{f^\vee}} $ and allowing $f$ to be a symbol
of H\"ormander class $\SmA$ we obtain a class of multipliers extending
the pseudodifferential multipliers \`a la Connes \cite{Con:CAG} and 
Baaj's \cite{Baa:CPDII,Baa:CPDI}. Later we will also need the 
so called classical ($1$--step polyhomogeneous) symbols
$f\in\CSmA$ which have an asymptotic expansion
\[
  f\sim \sum_{j=0}^\infty f_{m-j}
\]
with
$f_{m-j}(\gl\xi)=\gl^{m-j}\cdot f_{m-j}(\xi), |\xi|\ge 1, \gl\ge 1$. 

Thus for $f\in \SmA$ we obtain a well--defined multiplier $P_f$
acting on the pre-$C^*$--module $\SRnA$ with \emph{complete symbol}
$f$. The usual stationary phase arguments (e.g.
\cite[\S~I.3]{Shu:POS}) then allow to prove that the space
$\pdo^\bullet_\gs(\R^n,\Ainf)=\bigcup_{n\in\Z} 
\pdo^m_\gs(\R^n,\Ainf)$ of twisted
pseudodifferential multipliers (as well as its classical counterpart
where the symbols $f$ are $1$--step polyhomogeneous) is a
$*$--algebra.

For symbols $f\in\SmA, g\in \SmAvar{m'}$
the composition $P_f\circ P_g$ is a pseudodifferential multiplier with
symbol $h\in\SmAvar{m+m'}$ and $h$ has the asymptotic expansion
\begin{equation}\label{EqPsiDOproduct}  
  h(t)\sim \sum_{\gamma} \frac{i^{-|\gamma|}}{\gamma!} (\pl^\gamma f)(t)
  \pl_y^\gamma\big\vert_{y=0} \Bl \ga_{-y}\bl g(t+By)\br\Br.
\end{equation}

Furthermore, $P_f^*$ is a pseudodifferential multiplier with symbol
\begin{equation}\label{EqSymAdjoint}  
 \sigma(P_f^*)\sim \sum_\grg \frac{1}{\grg!} \pl_t^\grg \delta^\grg f(t)^*.
\end{equation}
Here $\delta^\grg$ denotes the basic derivative on $\sA$ induced by
the flow $\ga$:  For $a\in\Ainf$ and a multiindex $\gamma\in\Z_+^n$
it is defined by 
\begin{equation}\label{EqDeltaOnA}
  \gd^\gamma  a := i^{|\gamma|} \pl_x^\gamma\bigl|_{x=0} \ga_x(a)=
   i^{-|\gamma|} \pl_x^\gamma\bigl|_{x=0} \ga_{-x}(a).
\end{equation}
$\gd^\gamma$ plays the role of the partial derivative
$i^{-|\gamma|}\frac{\pl^\gamma}{\pl x^\gamma}$. 

\subsection{Differential multipliers} 
In the standard calculus differential operators are characterized as
those pseudodifferential operators whose complete symbols are
polynomial in the covariables $\xi$. Adopting this as a definition for
\emph{differential multipliers} it turns out that the (multiplier
counterparts) of the natural first and second order differential
operators discussed in Section \ref{MLSec2} are differential
multipliers in this sense.

Somewhat more formally we call $P_f\in\pdo^\bullet_\gs(\R^n,\Ainf)$ 
a \emph{differential multiplier} of order $m$ if
\begin{equation}\label{EqDiffOpSym}
  f(\xi)= \sum_{|\gamma|\le m} a_\gamma \xi^\gamma; \quad
  a_\gamma\in\Ainf,
\end{equation}
$f\in \Ainf[\xi_1,\ldots,\xi_n]$ is a polynomial of degree at most
$m$. Here the sum runs over all multiindices $\gamma\in \Z_+^n$ with
$|\gamma|\le m$.  Clearly, polynomials in $\xi$ are $1$--step
polyhomogeneous and hence differential multipliers are
\emph{classical} pseudodifferential multipliers.

Recall that in the ordinary pseudodifferential calculus the symbol of the basic
derivatives $i^{-|\gamma|} \pl_x^\gamma$ is given by $\xi^\gamma$.
Therefore, for a multiindex $\gamma$ we put
$\ud^\grg:=P_{\xi^\gamma}$. Explicitly, we find from \Eqref{EqPsiDOb}
for $u\in \SRnA$
\begin{equation} 
\begin{split}
(\ud^\gamma u)(x)  & = (P_{\xi^\grg}u)(x)= \int_{\R^n} e^{i\inn{x,y}}  \bl \xi+Bx\br^\grg \hat u(\xi) \, \dsl \xi \\
      & = i^{-|\grg|} \pl_y^\grg\bigl|_{y=0} 
            \int_{\R^n} e^{i \inn{x+y,\xi+Bx}} \hat u(\xi)\, \dsl\xi\\
      & = i^{-|\grg|} \pl_y^\grg\bigl|_{y=0}\bl e(x,y) u(x+y) \br   
       = i^{|\grg|} \pl_y^\grg\bigl|_{y=0} \uU_y u(x).
\end{split}
  \label{EqML4}
\end{equation}

It is important to note that due to the twisting in
general $\ud^{\grg}\ud^{\grg'}\not= \ud^{\grg+\grg'}$, as can
be seen either directly  or by the just proved product formula.

 As in the ordinary pseudodifferential calculus it is
in general not true that $P_f^*=P_{f^*}$. However,
$\sigma(P_f)^*=\sigma(P_f^*)\mod \SmAvar{m-1}.$

Furthermore, $\ud^\gamma$ is formally self-adjoint and thus for
any \emph{differential} multiplier we have indeed $P_f^*=P_{f^*}$.

\subsection{Differential multipliers of order $1$ and $2$}
 \label{SSSDiffOp12}
We look more closely at the most relevant case of differential
multipliers of order $1$ and $2$.  Let $e_j, j=1,\ldots, n$ 
be the canonical basis vectors of $\R^n$. 
We abbreviate $\ud_j:=\ud^{e_j}$ and recall that $b_{jk}$
denotes the entries of the skey--symmetric structure matrix
of the twisting \Eqref{Eq1206041}. Then by \Eqref{EqML4}
\begin{equation}\label{eq:1209128}  
 \ud_j u(x)= i\ii\pl_{y_j}\bigl|_{y=0} e^{i\inn{Bx,y}} u(x+y)
           = \bigl(\frac 1i \pl_{x_j} + b_{jl} x_l\bigr) u(x),
\end{equation}
were summing over repeated indices is understood. Thus
\begin{align}\label{eq:1209129}  
 \ud_j \ud_k &= -\pl_{x_j}\pl_{x_k} - i b_{js} x_s \pl_{x_k} - i  b_{ks} x_s
  \pl_{x_j} - i b_{kj} + b_{js} b_{kr} x_s x_r.
\end{align}
In particular we have the ``curvature identity''
\begin{equation}\label{EqMLCurvature}
[\ud_j,\ud_k] = 2 i \, b_{jk}.
\end{equation}

The twisting and the non-commutativity has an intersting effect on
the symbol calculus. The symbol of $\ud_j\cdot \ud_k$ is not
$\xi_j\cdot \xi_k$ but rather it is a consequence of the
formula \Eqref{EqPsiDOproduct} that
\begin{equation}\label{EqSymProd}   
     \sigma\bigl(\ud_j \cdot \ud_k \bigr)  
       = \xi_j \cdot \xi_k + i b_{jk} 
       = \sigma\bigl(\ud^{e_j+e_k}\bigr)+ i b_{jk},
\end{equation}
hence $ \ud^{e_j+e_k} = \ud_j\cdot \ud_k - i b_{jk}$. From this
the curvature identity \Eqref{EqMLCurvature} also follows.

\subsubsection{Differential multipliers in dimension $n=2$} 
\label{SSSDiffOp2}
Specializing further to dimension $n=2$ it is most convenient to
make use of the complex Wirtinger derivatives. Furthermore,
the structure matrix $b_{jk}$ has only one interesting entry $b_{12}$.
Fixing $\tau\in\mathbb{C}$ with $\Im\tau>0$ (a complex structure!) 
we have the following basic differential multipliers:
\[ 
\begin{split}
 \udtau &:= \ud_1 + \taubar \ud_2,  
     \quad \udtau^*=\ud_1 + \tau \ud_2,\quad
                     \ud_1 := \ud^{1,0}, \ud_2:= \ud^{0,1} \\
    [\udtau,\udtau^*] & = - 4\cdot \Im\tau\cdot b_{12}=:c_\tau,\\
          \uDelta_\tau & := \frac 12 (\udtau^*\udtau+\udtau\udtau^*)
        = \ud_1^2+ |\tau|^2 \ud_2^2 + \Re\tau ( \ud_1 \ud_2+\ud_2\ud_1).
\end{split}        
\]
We will first analyze these operators acting as multipliers on the Hilbert
module completion of $\SRnA$. Lateron we will have to pass to their concrete
counterparts acting on the Heisenberg modules.

\section{The resolvent expansion and trace formula} \label{ResExp}

The resolvent trace, or equivalently the heat trace, expansion for
second order Laplace type operators goes back at least to
Minakshisundaram and Pleijel \cite{MinPle:SPE}. Via Karamata's
tauberian theorem there is a connection to the eigenvalue counting
function whose asymptotic analysis is quite subtle. The best remainder
term for the counting function of general elliptic operators led
H\"ormander to develop his beautiful theory of Fourier integral 
operators \cite{Hor:FIOI}. Later the resolvent trace (aka heat
equation) method led to the development of local index theory
\cite{ABP73} with an enormous flow of publications.

In our opinion, by now the most streamlined approach to the
resolvent expansion of elliptic \emph{differential} operators
is the calculus of parameter dependent pseudodifferential
operators which essentially goes back to Seeley's seminal
complex powers paper \cite{See:CPE} and which is presented
very nicely in Shubin's book \cite[\S~II.9]{Shu:POS}.
\footnote{In Seeley's paper a subtle oversight caused a certain
confusion which, at least among non-experts, seems to exist
to this day. The resolvents of elliptic pseudodifferential
operators in general only belong to a ``weakly parametric``
calculus. This difference between the resolvent calculi
for differential resp. true pseudodifferential operators was
  clarified almost 30 years after Seeley's original paper
  \cite{GruSee:WPP}.}
We will come back to this soon. Our goal here is to show
that this calculus carries over to twisted pseudodifferential
multipliers and that the second coefficient in the expansion
can be calculated quite easily without any computer aid.
 
We consider the differential  multiplier
$P= P_{\eps_1 , \eps_2} := k^2 \uDelta_\tau + \eps_1 (\dtau k^2) \udtau^* + \eps_2
  (\dtau^*k^2) \udtau + a_0$,
where $a_0\in\sA$ and $\eps_1, \eps_2$ are real parameters. 
This multiplier contains all conformal Laplace type multipliers,
which occur on Heisenberg modules over noncommutative tori, as special
cases. The symbol of $P$ takes the form
$\sigma_P(\xi):= a_2(\xi) + a_1(\xi) +a_0$,
where $a_0\in\Ainf$ is the same as above and
\begin{align*}
    a_2(\xi) &= k^2 |\xi_1+\taubar \xi_2|^2 =: k^2 |\eta|^2, \\
    a_1(\xi) &= \eps_1 (\dtau k^2) \etabar + \eps_2 (\dtau^* k^2) \eta,
             &\eta &:= \xi_1+\taubar \xi_2, \\
             &=: \varrho_1 \etabar + \varrho_2 \eta, &\vrho_1 &:= \eps_1
  \dtau k^2, \varrho_2 := \eps_2 \dtau^* k^2.
\end{align*}

The resolvent $(P-\gl)\ii$ belongs to the parameter dependent
pseudodifferential calculus and therefore its symbol has a
polyhomogeneous expansion 
$\sigma_{(P-\gl)\ii}\sim b_{-2}+ b_{-3}+b_{-4}+\ldots$,
where $b_{-k}(\xi,\gl)\in\Ainf$ depends smoothly on $(\xi,\gl)$
and is \emph{homogeneous} of degree 
$-k$: $b_{-k}(r\xi, r^2\gl) = r^{-k} b(\xi,\gl)$. 
As a consequence we obtain for the $a\in\Ainf$ with respect
the \emph{dual trace} \Eqref{EqDualTrace} $\hat\varphi_0$
($\varphi_0$ is the invariant trace on $\Ainf$)
an asymptotic expansion
\begin{equation}
  \varphi_0 \Bl e^{-tP}\Br\sim_{t\searrow 0} \sum_{j=0}^\infty
     a_{2j}(P,a) t^{j-1},
\end{equation}
where it follows from the homogeneneity
\footnote{Note that heat/resolvent invariants are enumerated from
$0$. We are after $a_2$ which is the second nontrivial heat invariant, as
$a_1$ is always $0$ for differential operators, but in the counting of the
recursion system it is the third term.} 
\begin{equation}
  \begin{split}
    a_{2j}(P,a) &= \int_{R^2} \int_C e^{-t\gl} \varphi_0\bl
    (b_{-2j}(\xi,gl)\br \dsl \gl \dsl \xi\\
      & = \int_{\R^2} \varphi_0\Bl b_{-2j-2}(\xi,-1)\Br \dsl\xi
          \int_C e^{-t\gl} (-\gl)^{-j} \dsl\gl.
  \end{split}
\end{equation}
Here $C$ is a contour in the complex plane encircling the positive
semiaxis clockwise such that 
$\int_C e^{-t\gl} (r-\gl)\ii \dsl\gl = e^{-tr}$.
The second line is a consequence of the homogeneity of the $b_{-k}$
(see \cite[\S 6]{ConMos:LIF}). For the second nontrivial heat 
coefficient
one therefore obtains up to a sign (see loc. cit.)
\[
    a_{2}(P,a) = \int_{\R^2} \varphi_0( b_{-4}(\xi,-1) ) \dsl\xi.
\]
Due to this formula, it will be convenient to compute $b_{-4}(\xi,-1)$
modulo functions of total $\xi$--integral $0$. Up to a
function of total $\xi$--integral $0$ we have the following closed formulas
for the first three terms in the symbol expansion of $(P-\gl)\ii$:
\begin{align*}
 b_{-2}  & = b  = (k^2 |\eta|^2-\gl)\ii,
   \quad b_{-3}    = - b k^2 \bl \eta \dtau^*+ \etabar \dtau\br b - b a_1 b,
                  \\
  b_{-4}    & = \bl 2b k^2 |\eta|^2 - 1 - \eps_1-\eps_2\br b k^2
                    \Lapl_\tau b 
    + \gl b k^2 \bl (\dtau^* b) ( \dtau b) 
                    + (\dtau b)(\dtau^* b)\br 
                  \\
    &\quad + \eps_1\cdot \gl b (\dtau k^2) b \dtau^* b         
              + \eps_2 \cdot \gl b (\dtau^* k^2) b \dtau b   
                 \\
           &\quad + \eps_1\eps_2\cdot 
               |\eta|^2 b\cdot \bl (\dtau k^2) b (\dtau^* k^2)
               +(\dtau^* k^2) b \dtau k^2 \br \cdot b - b a_0 b.
\end{align*}

The proof is straightforward, completely computer free, and fits
on two pages, cf. \cite[\S~3.3]{LM2015}.

\subsection{Second heat coefficient} \label{SecHC}

Integrating $b_{-4}$ over $\xi$ is still a little involved
and it requires the Rearrangement Lemma \cite[\S~6.2]{ConMos:LIF}.
This was recast and generalized in \cite{Les2014}. The calculus
of divided differences allows to compute the many explicit
integrals in a systematic way. As a result 
there exist entire functions $K(s), H^\Re(s,t), H^\Im(s,t)$, 
such that with $h:=\log k^2$  the second heat coefficient of $P$ 
(w.r.t. the natural dual trace on the twisted crossed product) takes the form
\begin{multline*}
  a_2(P,a) = \frac{1}{4\pi |\Im\tau|} \varphi_0\Biggl[ a 
         \Bl K(\nabla)(\Lapl_\tau h)  - k^{-2} a_0 \\
  + H^\Re(\nabla^{(1)},\nabla^{(2)})\bl \square^\Re(h) \br
   + H^\Im(\nabla^{(1)},\nabla^{(2)})\bl \square^\Im(h) \br \Br
            \Biggr].
\end{multline*}
Here, $\square^{\Re/\Im}(h):= \frac 12 \bl
\dtau h\cdot \dtau^* \pm  \dtau^* h\cdot \dtau h\br$,
$\nabla = -\operatorname{ad}(h)$, and $\nabla^{(i)}$ signifies
that it acts on the $i$-th factor (\textit{cf}. \cite{CM2014}, 
\cite{Les2014}).

{The functions} $K, H^\Re, H^\Im$ depend only on $P$ but not on
$\tau$. They can naturally be expressed in terms of simple divided
divided differences of $\log$.

\subsection{Effective pseudodifferential operators and trace formulas} 
\label{sec-Epo}
%
We consider the noncommutative torus $\sA_\theta$ with generators
$U_1, U_2$ and normalized trace $\varphi_0$. 
Let $f:\R^2\to \Ainf$ be a symbol function (or Schwartz function) 
of sufficiently low order. Recall the trace \Eqref{EqDualTrace} of
the \emph{multiplier} $P_f$:
$\hat\varphi(P_f) = \int_{\R^2} \varphi_0 (f(x) ) dx$. However,
the multiplier $P_f$ is canonically represented as an operator
on the GNS space $L^2(A_\theta,\varphi_0)$ by $\Op(f)=\int_{\R^2}
f^{\vee}(x) \pi_x
dx$, where $\pi_x(U_1^{n_1}U_2^{n_2} )= e^{ i \inn{x,n} } U_1^{n_1}
U_2^{n_2}$ is the unitary which implements the natural $\R^2$--action
on $A_\theta$, cf. \Eqref{actT}. $\Op(f)$ acts as a trace class
operator on $L^2(A_\theta,\varphi_0)$. More concretely,
one computes $\Op(f) U_1^{n_1}U_{n_2} = f(-n_1,-n_2) U_1^{n_1} U_2^{n_2}$.
Since $(U_1^{n_1}
U_2^{n_2})_{n\in\Z^2}$ is an orthonormal basis of
$L^2(\Ainf,\varphi_0)$ we obtain the trace formula
\begin{equation} \label{EqMLTraceFormula}
  \begin{split}
    \Tr(\Op(f) )  & = \sum_{n\in\Z^2} \inn{ U_1^{n_1}U_2^{n_2} ,
    f(-n_1,-n_2)U_1^{n_1}U_2^{n_2} }\\
      & = \sum_{n\in\Z^2} \varphi_0\Bl \bl U_1^{n_1}U_2^{n_2}\br^*
                f(-n_1,-n_2) U_1^{n_1} U_2^{n_2} \Br\\
      & = \sum_{n\in\Z^2} \varphi_0 \bl f(n) \br.
  \end{split}
\end{equation}
This looks a little different from the formula for the dual trace.
However, for a parameter dependent symbol $f(x,\gl)$ we can
take advantage of the Poisson summation formula. Then we find
\begin{equation}
  \begin{split}
    \sum_{n\in\Z^2} f(k,\gl) &= \sum_{k\in\Z} \hat f(2\pi k,\gl) \\
        & = \hat f(0,\gl) + \sum_{k\in\Z^2\setminus\{0\}} \hat f(2\pi
        k,\gl) \\
    & = \int_{\R^2} f(\xi,\gl) d\xi + O(\gl^{-N})
  \end{split}
\end{equation}
for any $N$. The latter follows from integration by parts in the
  Fourier transform and the symbol estimates.

Thus the upshot is that for trace class symbols in the parameter
dependent calculus the multiplier trace and the trace in the Hilbert
space representation coincide only \emph{asymptotically}. However,
for computing heat and resolvent trace asymptotics this is good
  enough. 

Furthermore, this observation has a far reaching generalization.  
Namely, the effective implementation of the pseudodifferential calculus amounts to
passing from its realization on multipliers to a direct action on projective
representation spaces (Heisenberg modules) or on $L^2(\sA,\gvf_0)$ itself.
More concretely, let  $\spi:G\to \cL(\cH)$ be a projective unitary
representation of $G=\R^n\times(\R^n)^\wedge$. For a symbol $f\in \SmA$ the
assignment $\SmA\ni f\mapsto \Op(f):= \int_G f^\vee(y) \spi(y) dy$ represents
pseudodifferential multipliers as concrete operators in $\cH$.

By exploiting the representation theory of the Heisenberg group we are able to
relate the Hilbert space trace of parameter dependent pseudodifferential
operators to the trace of the corresponding multiplier acting on $\SRnA$.  
For details see \cite[\S~5 and Appendix A]{LM2015}.



 



 

\bibliography{mlbib.bib,localbib.bib}
\bibliographystyle{amsalpha-lmp}
\end{document}